\documentclass[10pt,a4paper,titlepage]{article}
\usepackage[latin1]{inputenc}
\usepackage{amsmath}
\usepackage{amsfonts}
\usepackage{amssymb}
\usepackage{amsthm}
\usepackage[small,nohug,heads=LaTeX]{diagrams}
\usepackage{amsfonts}
\usepackage{amssymb}
\usepackage[dvips]{graphicx}
\usepackage{fancyhdr}
\usepackage{float}
\diagramstyle[labelstyle=\scriptstyle]
\makeatletter
\def\imod#1{\allowbreak\mkern10mu({\operator@font mod}\,\,#1)}
\makeatother
\numberwithin{figure}{section}
\theoremstyle{plain}
\newtheorem{thm}{Theorem}[section]
\newtheorem*{prop*}{Proposition}
\newtheorem*{thm*}{Theorem}
\newtheorem*{thmA}{Theorem \ref{thm: A}}
\newtheorem*{thmB}{Theorem \ref{thm: B}}
\newtheorem*{thmC}{Theorem \ref{thm: C}}
\newtheorem{prop}[thm]{Proposition}
\newtheorem{lem}[thm]{Lemma}
\newtheorem{cor}[thm]{Corollary}
\theoremstyle{definition}

\newtheorem{dfn}[thm]{Definition}
\newtheorem*{dfn*}{Definition}
\newtheorem{notation}[thm]{Notation}
\theoremstyle{remark}
\newtheorem{rmk}[thm]{Remark}

\def\C{\mathbb{C}}

\def\N{\mathbb{N}}
\def\Z{\mathbb{Z}}

\def\K{\mathbb{K}}
\def\S{\mathbb{S}}
\def\emb{\hookrightarrow}
\def\Wlog{Without loss of generality }

\def\rep{representation }
\def\reps{representations }
\def\Torelli{\mathrm{\overline{IA}}}
\def\iff{if and only if }

\def\GL{\mathrm{GL}}

\def\into{\hookrightarrow}

\newcommand{\Out}[1]{\mathrm{Out}(F_{#1})}
\newcommand{\Aut}[1]{\mathrm{Aut}(F_{#1})}
\newcommand{\Inn}[1]{\mathrm{Inn}(F_{#1})}

\begin{document}
\textsc{\begin{LARGE}\begin{center} Outer automorphism groups of free groups: linear and free representations
\end{center}\end{LARGE}}

\medskip

\begin{center}
Dawid Kielak\footnotemark

\today
\end{center}

\medskip

\begin{center}
\begin{minipage}{0.65\textwidth}
\textsc{Abstract.} We study homomorphisms between $\Out{n}$ and $\Out{m}$ for $n \geqslant 6$ and $m < {n \choose 2}$, and conclude that if $m \neq n$ then each such homomorphism factors through the finite group of order 2. In the course of the argument linear \reps of $\Out{n}$ in dimension less than ${{n+1} \choose 2}$ over fields of characteristic zero are completely classified. It is shown that each such \rep has to factor through the natural projection $\Out{n} \to \GL_n(\Z)$ coming from the action of $\Out{n}$ on the abelianisation of $F_n$. We obtain similar results about linear representation theory of $\Out 4$ and $\Out 5$.
\end{minipage}
\end{center}

\bigskip

\footnotetext[1]{The author was supported by the EPSRC of the United Kingdom}


\section{Introduction}
In this paper we study the problem of existence of homomorphisms between outer automorphism groups of finitely generated free groups of different rank.
In the free abelian case the automorphism groups (and in fact the outer automorphism groups) are of the form $\GL_n(\Z)$. Among many other properties, these groups easily embed into each other: to be precise, if $m \geqslant n$, we can construct an embedding $\GL_n(\Z) \emb \GL_m(\Z)$ just by mapping every matrix to the $n \times n$ upper-left corner, and then completing the matrix by putting an identity in the lower-right corner. Actually, the more general problem of understanding homomorphisms $\GL_n(\Z) \to \GL_m(\Z)$ is completely solved (mostly thanks to Margulis's superrigidity).

Similarly to the free abelian case, we can construct embeddings $\Aut{n} \emb \Aut{m}$ between automorphism groups of finitely generated free groups (with $m \geqslant n$ as before), by choosing a free factor in $F_m$ isomorphic to $F_n$.

The situation becomes far less obvious when we focus on the outer automorphism groups of finitely generated free groups. Until recently very little was known about possible embeddings between groups of the form $\Out{n}$. The positive results known to the author were obtained by Aramayona, Leininger and Souto~\cite{aramayonaetal2009}, Bogopol'skii and Puga~\cite{bogopol'skiipuga2002} (of which a slightly stronger version was proven by Bridson and Vogtmann~\cite{bridsonvogtmann2011}) and Khramtsov~\cite{khramtsov1990}.

Some negative results were also obtained:
\begin{enumerate}
\item Khramtsov~\cite{khramtsov1990} has proven that $\Out{n}$ never embeds into $\Out{n+1}$ (if $n>1$);
%
%
\item Bridson and Vogtmann~\cite{bridsonvogtmann2003} have shown that for any $n\geqslant 3$ and $m<n$, there exist no embeddings $\Out{n} \emb \Out{m}$. In fact their result shows that the image of any homomorphism $\Out{n} \to \Out{m}$ is contained in a copy of $\Z_2$, the cyclic group of order 2;

\item Bridson and Vogtmann~\cite{bridsonvogtmann2011} have also shown that the image of any homomorphism $\Out{n} \to \Out{m}$ is of size at most 2 whenever $n$ is at least $9$, $m \neq n$, and $m \leqslant 2n-2$ when $n$ is odd, or $m \leqslant 2n$ when $n$ is even.
\end{enumerate}

In particular, the last result gives an answer to a question of Bogopol'skii and Puga, who in~\cite{bogopol'skiipuga2002} conjectured that there always exist embeddings $\Out{n} \emb \Out{2n}$. Bridson and Vogtmann have shown that there are no embeddings of this form provided that $n \geqslant 9$ and $n$ is even. We extend their result to the case $n \geqslant 6$, independently of the parity of $n$. We obtain

\begin{thmA}
Let $n, m \in \N$ be distinct, $n \geqslant 6$, $m< {n\choose2}$, and let $\phi \colon \Out{n} \to \Out{m}$ be a homomorphism. Then the image of $\phi$ is contained in a copy of $\Z_2$, the finite group of order two.
\end{thmA}

A result of Bridson and Farb~\cite{bridsonfarb2001} allows us to extend this result further. We prove

\begin{thmC}
Let $n, m \in \N$ be distinct, with $n$ even and at least $6$. Let $\phi \colon \Out{n} \to \Out{m}$ be a homomorphism. Then the image of $\phi$ is finite, provided that
\[ {n \choose 2} \leqslant m< {{n+1}\choose2}. \]
\end{thmC}

The question of finding $n$ for which we have $\Out n \into \Out {2n}$ has not been fully answered. The case $n=1$ is trivial, and Khramtsov~\cite{khramtsov1990} has shown that there exists an embedding $\Out 2 \into \Out 4$. Our result (and the work of Bridson--Vogtmann~\cite{bridsonvogtmann2011}) dealt with the case $n \geqslant 6$. The cases $n \in \{3,4,5\}$ remain unanswered, however the author has approached the solution to the $n=3$ case in \cite{kielak2011(2)}, where he shows that $\Out 3$ does not embed into $\Out 5$.

The general strategy of this paper consists of two steps (see~\cite{bridsonvogtmann2011} for a similar approach). Firstly, we investigate the low-dimensional representation theory of $\Out{n}$, which in particular enables us to prove

\begin{thmB}
Let $\K$ be a field of characteristic equal to zero or greater than $n+1$.
Suppose $\phi  \colon \Out{n} \to \GL(V)$ is an $m$-dimensional $\K$-linear \rep of $\Out{n}$, where $n \geqslant 6$ and $m< {{n+1} \choose 2}$. Then $\phi$ factors through the natural projection $p  \colon \Out{n} \to \GL_n(\Z)$.
\end{thmB}

We then concentrate on obtaining information about allowed representations of a carefully chosen finite subgroup of $\Out{n}$. We use a result proven independently by Culler, Khramtsov and Zimmermann, to realise the action of our finite group on the conjugacy classes of $F_m$ as induced by an action on a finite graph. Comparing the representation theory with the action on the homology of this graph will yield the result.

At this point we owe an explanation to the reader. The use of torsion is indeed crucial -- it has been shown by Bridson and Vogtmann ~\cite[Proposition 2.5]{bridsonvogtmann2011} that for each positive $n$ there exists a finite index subgroup $\Gamma < \Out{n}$ which embeds into $\Out{m}$ for $m = 2n-1$.

The interplay between linear representations of $\Out{n}$ and homomorphisms $\Out{n} \to \Out{m}$ prompted us to coin the term `free representations' to describe homomorphisms $G \to \Out{m}$ for any group $G$.

\textsl{Acknowledgements.}
The author wishes to express deep gratitude to his supervisor, Martin Bridson, for the help obtained during the preparation of this paper. We also wish to thank the referee for his comments on exposition and David Craven for helpful conversations and pointing out Theorem~\ref{mapfromglnz2} to us.

\section{Notation and preliminaries}
\label{secnotation}
Let us first establish some conventions and definitions:

\begin{dfn}
\label{graphdef}
We say that $X$ is a \emph{graph} \iff it is a 1-dimensional CW complex. The 1-cells of $X$ will be called \emph{edges}, the 0-cells will be called \emph{vertices}. The sets of vertices and edges of a graph will be denoted by $E(X)$ and $V(X)$ respectively. The points of intersection of an edge with the vertex set are referred to as \emph{endpoints} of the edge.

We will equip $X$ with the standard path metric in which the length of each edge is 1.

Given two graphs $X$ and $Y$, a function $f  \colon X \to Y$ is a \emph{morphism of graphs}  \iff $f$ is a continuous map sending $V(X)$ to $V(Y)$, and sending each open edge in $X$ either to a vertex in $Y$ or isometrically onto an open edge in $Y$.

When we say that a group $G$ \emph{acts on a graph} $X$, we mean that it acts by graph morphisms.

We say that a graph $X$ is \emph{directed} \iff it comes equipped with a map $o  \colon E(X) \to X$ such that $o(e)$ is a point on the interior of $e$ of distance $\frac13$ from one of its endpoints. We also define $\iota, \tau  \colon E(X) \to V(X)$ by setting $\tau(e)$ to be the endpoint of $e$ closest to $o(e)$, and $\iota(e)$ to be the endpoint of $e$ farthest from $o(e)$. Note that we allow $\iota(e) = \tau(e)$.

%
The \emph{rank} of a connected graph is defined to be the size of a minimal generating set of its fundamental group (which is a free group). 
\end{dfn}

\begin{rmk}
Unless specified otherwise, the graphs we will be dealing with will be connected and non-trivial (that is with at least one edge).
%
%
\end{rmk}

\begin{notation}
Let $G$ be a group. We will adopt the following notation:
\begin{itemize}
\item $1$ denotes the identity element of $G$;
\item for two elements $g,h \in G$, we define $g^h=h^{-1}gh$;
\item for two elements $g,h \in G$, we define $[g,h]=ghg^{-1}h^{-1}$;
\item $G \curvearrowright X$ denotes a left action of $G$ on a set $X$, and $g.x$ is the image of $x \in X$ under $g \in G$;
\item the finite cyclic group of order $k$ will be denoted by $\Z_k$;
\item the free group of rank $n$ will be denoted by $F_n$.
\end{itemize}
\end{notation}

\begin{dfn}
Let us introduce the following notation for elements of $\Aut{n}$, the automorphism group of $F_n$, where $F_n$ is the free group on the set $\{ a_1, \ldots, a_n\}$:
\begin{eqnarray*}
\epsilon_i  \colon \left\{ \begin{array}{cccl} a_i & \mapsto & a_i^{-1}, \\ a_j & \mapsto & a_j, & j \neq i \end{array} \right., & &
\sigma_{ij}  \colon \left\{ \begin{array}{cccl} a_i & \mapsto & a_j, \\ a_j & \mapsto & a_i, & \\ a_k & \mapsto & a_k, & k \not \in \{i,j\} \end{array} \right., \\
\rho_{ij}  \colon \left\{ \begin{array}{cccl} a_i & \mapsto & a_i a_j, \\ a_k & \mapsto & a_k, &  k \neq i \end{array} \right., & &
\lambda_{ij}  \colon \left\{ \begin{array}{cccl} a_i & \mapsto & a_j a_i, \\ a_k & \mapsto & a_k, &  k \neq i \end{array} \right..
\end{eqnarray*}
Let us also define $\Delta = \prod_{i=1}^n \epsilon_i$ and
\begin{displaymath}
\sigma_{i (n+1)}  \colon \left\{ \begin{array}{cccl} a_i & \mapsto & a_i^{-1}, \\ a_j & \mapsto & a_j a_i^{-1}, & j \neq i \end{array} \right..
\end{displaymath}
\end{dfn}

We are going to use the same symbols to denote the images of those elements under the natural projection $\Aut n \to \Out n$.

Below we give an explicit presentation of $\Out{n}$, the outer automorphism group of $F_n$. It is the Gersten's presentation (see~\cite{gersten1984}) adapted to our notation and conventions (conjugation, commutators, action on the left etc.). Compare~\cite{bridsonvogtmann2011}.

\begin{prop}[(Gersten's presentation)]
\label{gerstenpres}
Suppose $n \geqslant 3$. The group $\Out{n}$ is generated by $\{ \epsilon_1, \rho_{ij}, \lambda_{ij} \mid i,j = 1,\ldots,n , \ i \neq j \}$, with relations
\begin{itemize}
\item $[\rho_{ij},\rho_{kl}]=[\lambda_{ij},\lambda_{kl}]=1$ for $k \not \in \{i,j\}, l \neq i$;
\item $[\lambda_{ij},\rho_{kl}]=1$ for $k \neq j, l \neq i$;
\item $[\rho_{ij}^{-1},\rho_{jk}^{-1}]=[\rho_{ij},\lambda_{jk}]=[\rho_{ij}^{-1},\rho_{jk}]^{-1}=[\rho_{ij},\lambda_{jk}^{-1}]^{-1}=\rho_{ik}^{-1}$ for $k \not \in \{i,j\}$;
\item $[\lambda_{ij}^{-1},\lambda_{jk}^{-1}]=[\lambda_{ij},\rho_{jk}]=[\lambda_{ij}^{-1},\lambda_{jk}]^{-1}=[\lambda_{ij},\rho_{jk}^{-1}]^{-1}=\lambda_{ik}^{-1}$ for $k \not \in \{i,j\}$;
\item $\rho_{ij} \rho_{ji}^{-1} \lambda_{ij}=\lambda_{ij} \lambda_{ji}^{-1} \rho_{ij}, \ (\rho_{ij} \rho_{ji}^{-1} \lambda_{ij})^4=1$;
\item $[\epsilon_1, \rho_{ij}] = [\epsilon_1, \lambda_{ij}]=1$ for $i,j \neq 1$;
\item $\rho_{12}^{\epsilon_1}=\lambda_{12}^{-1}, \ \rho_{21}^{\epsilon_1}=\rho_{21}^{-1}$;
\item $\epsilon_1^2=1$;
\item $\prod_{i\neq j} \rho_{ij} \lambda_{ij}^{-1} = 1$ for each fixed $j$.
\end{itemize}
\end{prop}
%

Note the action of $\Aut{n}$ on $F_n$ and $\Out{n}$ on the conjugacy classes of $F_n$ is \textbf{on the left}.

\begin{rmk}
Let us note that the following relations hold both in $\Aut{n}$ and $\Out{n}$:
$\rho_{ij}^{\epsilon_j} = \rho_{ij}^{-1}$ and
$\rho_{ij}^{\epsilon_i} = \lambda_{ij}^{-1}$, for all $i \neq j$.

It is now clear from the presentation that \[ \Out{n}= \big\langle \big\{ \rho_{ij}, \epsilon_i \mid i,j = 1,\ldots\,n, \ i\neq j \big\} \big\rangle. \]
\end{rmk}

%
%
\begin{dfn}
Let us define some finite subgroups of $\Out{n}$:
\begin{eqnarray*}
S_n &\cong&\langle \{\sigma_{ij} \mid i,j = 1,\ldots,n, \ i \neq j\} \rangle,\\
S_{n+1} &\cong&\langle \{\sigma_{ij} \mid i,j = 1,\ldots,n+1, \ i \neq j \} \rangle,\\
\Z_2^n \rtimes S_n\cong W_n& =& \langle \{ \epsilon_1, \sigma_{ij} \mid i,j=1\ldots,n, \ i \neq j \} \rangle, \\
\Z_2 \times S_{n+1}\cong G_n & =& \langle \{ \Delta, \sigma_{ij} \mid i,j=1\ldots,n+1, \ i\neq j \} \rangle.
\end{eqnarray*}
\end{dfn}

We do not give distinctive names to to the first two groups; instead, we will usually refer to them as respectively $S_n < W_n$ and $S_{n+1}<G_n$. More generally, whenever we mention $S_n$ or $S_{n+1}$ as subgroups of $\Out{n}$, we mean these two groups.

Note that we abuse notation by also using $S_n$ to denote the abstract symmetric group of degree $n$. We will denote its maximal alternating subgroup by $A_n$.

We will often talk about the natural action of $S_n$ and $A_n$ on $\{1,2,\ldots, n\}$. When doing so in the case of $S_n$, we will always mean the action in which
\[ \sigma_{ij}(k) = \left\{ \begin{array}{ccc} j & \textrm{ if } & k=i \\ i & \textrm{ if } & k=j \\ k & \textrm{ if } & k \neq i,j \end{array} \right. . \]
In the case of $A_n$, we will mean the restriction of the described action to $A_n < S_n$.

\begin{figure}
\begin{center}
\includegraphics[scale=0.9]{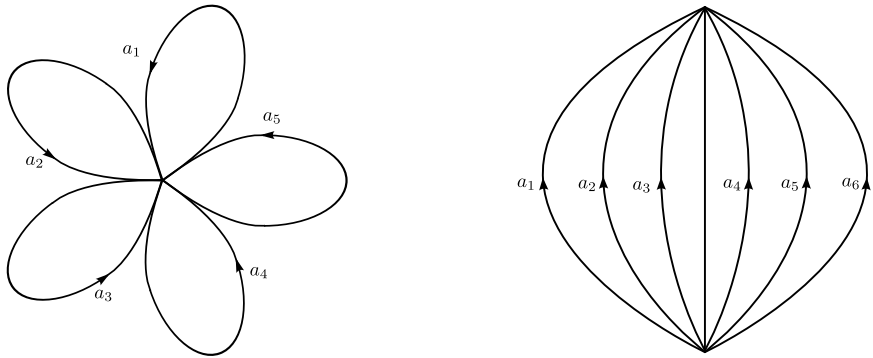}
\caption{The $5$-rose and $7$-cage graphs}
\label{rose and cage}
\end{center}
\end{figure}

Observe that the subgroup $W_n$ is the automorphism group of an $n$-rose, that is a graph with one vertex and $n$ edges, whereas the subgroup $G_n$ is the automorphism group of an $(n+1)$-cage, that is a graph with two vertices and ${n+1}$ edges, such that each edge has both vertices as its endpoints (see Figure~\ref{rose and cage}). Fixing an appropriate isomorphism between the fundamental groups of these graphs and $F_n$ induces the embeddings
$W_n, G_n < \Out{n}$.

Note that, if $i,j \leqslant n$, we have
\begin{displaymath}
\epsilon_i \sigma_{ij} = \lambda_{ij} \lambda_{ji}^{-1} \rho_{ij},
\end{displaymath}
and the subgroup ${S_n < \Out{n}}$ defined above acts on the sets
\begin{align*}
\{ \epsilon_i & \mid i=1\ldots,n\}, \\ \{ \rho_{ij} & \mid i,j=1\ldots,n, \ i \neq j\}, \textrm{ and} \\ \{ \lambda_{ij} & \mid i,j=1\ldots,n, \ i \neq j\}
\end{align*}
by permuting the indices in the natural way.

%
%
%
%
Let us also mention one useful fact (which is a consequence of a theorem by Mennicke~\cite{mennicke1965}):

\begin{prop}
\label{press}
The group $\mathrm{GL}_n(\Z_2)$ satisfies
\begin{displaymath}
\mathrm{GL}_n(\Z_2)= \Out{n} / \langle \! \langle \epsilon_{i} \mid i=1,\ldots,n \rangle \! \rangle.
\end{displaymath}
\end{prop}

\section{Linear representations of $\Out{n}$}

\begin{dfn}
Let us recall some basic terminology of representation theory of symmetric groups. Let $S_n$ be the symmetric group of degree $n$, and let $V$ be a $\K$-linear representation of the group, where $\K$ is a field of characteristic either 0 or greater than $n$. We will adopt the following terminology:
\begin{itemize}
\item if $\dim (V) =1$ and $S_n$ acts trivially on $V$, we say that $V$ is \emph{trivial};
\item if $\dim (V) =1$ and the representation $S_n \to \GL_1(\K) = \K^*$ has image $\{1,-1\}$, we say that $V$ is \emph{the determinant representation};
\item tensoring $V$ with the determinant representation adds the adjective \emph{signed} to the representation's name;
\item if $\dim (V) =n-1$ and the representation is the one induced on $\K$-homology by $S_n$ acting on the $n$-cage by permuting edges in the natural way, then we say that $V$ is \emph{the standard representation};
\item if $\dim (V) =n$ and the representation is the one induced on $\K$-homology by $S_n$ acting on the $n$-rose by permuting edges in the natural way, then we say that $V$ is \emph{the permutation representation}.
\end{itemize}
\end{dfn}

\begin{rmk}
All of these \reps can also be described using the correspondence between \reps of $S_n$ and partitions of $n$ (see e.g. \cite{fultonharris1991}); using this notation the trivial \rep corresponds to $(n)$, the determinant to $(1^n)$, the standard to $(n-1,1)$ and the permutation \rep is a sum $(n-1,1) \oplus (n)$.

Note that all but the last \rep are irreducible (the permutation \rep is a sum of two irreducible representations).
\end{rmk}

\begin{rmk}
The trivial and standard \reps are also irreducible \reps of $A_n < S_n$, the alternating group (for $n\geqslant 4$). The determinant \rep is trivial in this case, so there is no distinction between standard and signed standard representations of $A_n$.
\end{rmk}

Let us mention the following basic fact of representation theory.
\begin{lem}[(Schur's Lemma)]
\label{quotientgraphs}
\label{Schur's lemma}
Let $V$ be an irreducible \rep of a group $G$. Then any $G$-equivariant linear transformation $\phi  \colon V \to W$, where $W$ is a $G$-module, is either trivial or an isomorphism onto its image.
\end{lem}

We can now prove a useful result about \reps of $\Out{n}$.

\begin{lem}
\label{1-dim reps}
There are exactly two non-isomorphic $\K$-linear \reps of $\Out{n}$ of dimension 1, where $\K$ is a field of characteristic other than 2.
\begin{proof}
Let $\phi  \colon \Out{n} \to \GL_1(\K) \cong \K^*$ be a one-dimensional representation. Then
\[ \phi(\rho_{ij}) = \phi( [\rho_{ik}^{-1},\rho_{kj}^{-1}])^{-1} = 1, \]
since $\K^*$ is commutative. The same holds for any $\lambda_{ij}$, and hence $\phi$ factors as
\[ \begin{diagram} \Out{n}  & \rTo & \K^* \\ \dTo & & \uTo \\ \Out{n}/\langle \! \langle \rho_{ij}, \lambda_{ij} \mid i \neq j \rangle \! \rangle  & \cong & \Z_2  . \end{diagram}  \]
There are exactly two such maps; we are going to refer to the non-trivial one as the determinant map, and the corresponding \reps of $\Out{n}$ will also be called the determinant representation.
\end{proof}
\end{lem}

\begin{dfn}
\label{reps of W_n}
Let $V$ be a representation of $W_n \cong \Z_2^n \rtimes S_n$. Let $N = \{1,\ldots,n\}$. Define
\begin{itemize}
\item for each $I \subseteq N$,  $E_I = \{ v \in V \mid \epsilon_j v=(-1)^{\chi_I(j)}v\}$, where $\chi_I  \colon N \to \{0,1\}$ is the characteristic function of $I$;
\item $V_i = \bigoplus_{\vert I \vert =i} E_I$.
\end{itemize}
\end{dfn}

We will slightly abuse notation, and sometimes omit brackets and write $E_1$ for $E_{\{1\}}$, etc.

\begin{rmk}
Note that if $V$ is a $\K$-linear representation, where $\K$ is a field of characteristic other than 2, then we can simultaneously diagonalise the commuting involutions $\epsilon_1, \epsilon_2, \ldots, \epsilon_n$, and hence conclude that
\[ V = \bigoplus_{i=0}^n V_i = \bigoplus_{I \subseteq N} E_I. \]
\end{rmk}

\begin{lem}
\label{dimensionofV_i}
Let $V$ be a representation of $W_n$. Then, with the notation as above, $\dim V_i$ is divisible by $n \choose i$.
\begin{proof}
The symmetric group $S_n < W_n$ acts on $\{ \epsilon_1, \ldots , \epsilon_n \}$ by permuting the indices in the natural way. Hence its action on $V_i$ will be transitive on ${\{E_I \mid i=\vert I \vert \}}$. Therefore each $E_I$, for a fixed size of $I$, has the same dimension. The result follows by counting in how many ways we can pick a subset of size $i$ in $N$.
\end{proof}
\end{lem}

\begin{lem}
\label{diamond diagram}
Let $V$ be a representation of $\Out{n}$. Then, with the notation as above,
\begin{displaymath}
\rho_{ij}E_I \leqslant \bigoplus_{I\triangle J \subseteq \{i,j\}} E_J,
\end{displaymath}
where $A \triangle B$ denotes the symmetric difference of two sets, $A$ and $B$.
\begin{proof}
Let $v \in E_I$, and let $k \not\in \{i,j\}$. Note that
\[ [\rho_{ij}, \epsilon_k]=1\]
and that $v$ is an eigenvector of $\epsilon_k$ with eigenvalue $\mu_k \in \{ 1,-1\}$. Hence $\rho_{ij}(v)$ is also an eigenvalue of $\epsilon_k$ with eigenvalue $\mu_k$. Therefore
\[ v \in \bigoplus_{I\triangle J \subseteq \{i,j\}} E_J \]
as the space on the right hand side is the intersection of all the $\mu_k$-eigenspaces of elements $\epsilon_k$ for $k \not\in \{i,j\}$.
\end{proof}
\end{lem}

We will need a result about representations of a particular finite group due to Landazuri and Seitz \cite{landazuriseitz1974}:
\begin{thm}[(Landazuri, Seitz \cite{landazuriseitz1974})]
\label{mapfromglnz2}
Suppose $m < 2^{n-1}-1$, and $n>3$. Then every homomorphism $\GL_n(\Z_2) = \mathrm{PSL}_n(\Z_2) \to \mathrm{PGL}_m(\K)$ is trivial, provided that characteristic of $\K$ is not 2.
\end{thm}

In what follows, let us fix a field $\K$ of characteristic either 0 or greater than $n+1$.

\begin{prop}
\label{nchoose2case}
Suppose $V$ is an $m$-dimensional $\K$-linear \rep of $\Out{n}$, where $m< n(n-2)$, such that, with the notation of Definition~\ref{reps of W_n},
\[ \forall i \not \in \{0,1,n-1,n\} \colon V_i=\{0\}.\]
Assume also that $n \geqslant 6$ or that $n \geqslant 4$ and $\dim V_1 + \dim V_{n-1} = n$. Then $V$ decomposes as an $\Out{n}$-module as
\begin{displaymath}
V = V_0 \oplus V_1 \oplus V_{n-1} \oplus V_n,
\end{displaymath}
where the action of $\Out{n}$ on $V_0$ is trivial, and on $V_n$ it is via the determinant map. Moreover, as modules of $S_{n+1} < \Out{n}$, $V_1$ is a sum of standard, and $V_{n-1}$ of signed standard representations.
\begin{proof}
We are going to proceed in a number of steps.

\noindent \textbf{Step 0:} Let us first prove that $V = V_0 \oplus V_1 \oplus V_{n-1} \oplus V_n$ as an $S_{n+1}$- and $W_n$-module.

When $n \geqslant 5$,
Lemma~\ref{diamond diagram} tells us that \[\rho_{ij} (V_0 \oplus V_1) \leqslant V_0 \oplus V_1 \oplus V_2 \oplus V_3 = V_0 \oplus V_1\] since $V_2 = V_3=\{0\}$ by assumption.

When $n=4$ then $V_0 \oplus V_1 = V_0$ (as $V_1 = \{0\}$) or $V_3 = \{0\}$. In either case we conclude
\[\rho_{ij} (V_0 \oplus V_1) \leqslant V_0 \oplus V_1 \]
since $V_2 = \{0\}$.

Also, each $\epsilon_i$ keeps $V_0 \oplus V_1$ invariant, and therefore so does the entire group $\Out{n}$. Similarly $\Out{n}$ keeps $V_{n-1} \oplus V_n$ invariant.

The group $S_{n+1}$ commutes with $\Delta$, and preserves $V_0 \oplus V_1$ and $V_{n-1} \oplus V_n$, therefore it preserves each $V_i$. By construction we also see that each $\epsilon_j$ preserves the decomposition $V = V_0 \oplus V_1 \oplus V_{n-1} \oplus V_n$, and hence so does $W_n$.

\noindent \textbf{Step 1:} We claim that as $S_{n+1}$-modules, $V_1$ is a sum of standard representations, and $V_{n-1}$ is a sum of signed standard representations.

Let us look more closely at the representations of $S_{n-1}$ on $E_1$ and $E_{N \smallsetminus \{1\} }$, where $S_{n-1}$ is the stabiliser of $1$ when $S_{n}$ acts
on the indices of $\{ \epsilon_1, \ldots , \epsilon_n \}$.
Note that $E_1$ and $E_{N \smallsetminus \{1\} }$ are $S_{n-1}$-invariant, since $S_{n-1}$ commutes with $\epsilon_1$. The dimension of each of these representations is less than $n-2$ (by Lemma~\ref{dimensionofV_i} and our assumption on $m$). If $n \geqslant 6$ then these have to be  sums of trivial and determinant representations (see e.g.~\cite{rasala1977}). If $n \in \{4,5\}$ then \[\dim E_1 \in \{0,1\} \textrm{ and } \dim E_{N\smallsetminus \{1\}} \in \{0,1\}\] by assumption on dimensions of $V_1$ and $V_{n-1}$. Hence, as $S_{n-1}$-representations, $E_1$ and $E_{n-1}$ are sums of trivial and determinant representations.

Fix a basis $\{ b_1, \ldots, b_k \}$ of $E_1$, so that each $\langle b_i \rangle$ is $S_{n-1}$-invariant. We see that for each $i$, $\langle \sigma (b_i) \mid \sigma \in S_n \rangle$ is an $n$-dimensional representation of $S_n$, which has to be either the permutation or the signed permutation representation (since we know how $S_n$ acts on spaces $E_1, \ldots , E_n$). We immediately conclude, using the branching rule, (see e.g. \cite[Theorem 9.3]{james1978}) that the representation of $S_{n+1}$ on $V_1$ and $V_{n-1}$ is a sum of standard and signed standard representations.

Again we will focus on the subspaces $E_1$ and $E_{N \smallsetminus \{ 1 \} }$. We will only discuss the $E_1$ case, since the other case is analogous. Note that Lemma~\ref{diamond diagram} gives us
\begin{displaymath}
\rho_{ij}E_I \leqslant \bigoplus_{I\triangle J \subseteq \{i,j\}} E_J.
\end{displaymath}
Hence in particular
\begin{displaymath}
\rho_{ij} E_1 \leqslant E_1
\end{displaymath}
for all $i,j \neq 1$, since $E_{1,i} = E_{1,j} = E_{1,i,j} =\{0\}$, as $V_2 = V_3 = \{0\}$. But each $\rho_{ij}$ is an isomorphism, hence it has to be an isomorphism on $E_1$. Now the actions of $\rho_{23}$ and $\rho_{34}$ on $E_1$ are conjugate by the action of $\sigma_{24}\sigma_{34}$, which is trivial on $E_1$. Hence $\rho_{24} = [\rho_{34}^{-1}, \rho_{23}^{-1}]$ acts trivially on $E_1$. The same is true for $\lambda_{24}$ and $\lambda_{42}$, and hence $\sigma_{24} \epsilon_4 = \lambda_{24} \lambda_{42}^{-1} \rho_{24}$ acts trivially on $E_1$. Therefore the representation of $S_{n+1}$ on $V_1$ is a sum of standard representations, whereas the representation on $V_{n-1}$ is a sum of signed standard representations, which proves the claim.

Note that we have also shown that $\rho_{ij}$ acts as identity on $E_k$ and $E_{N \smallsetminus \{k\} }$ for each $k \not \in \{i,j\}$. This fact will turn out to be very useful in the remaining part of the proof.

\noindent \textbf{Step 2:} We now claim that $V_1$ and $V_{n-1}$ are $\Out{n}$-invariant.

In fact, we will only prove this claim for $V_1$, the $V_{n-1}$ case being analogous.

Note that the action of $A_n$ on $V_1$ gives isomorphisms $\iota_{ij}  \colon E_i \cong E_j$ for each $i,j$.
Let us consider $W \leqslant V_1$, an irreducible \rep of $S_{n+1}$. We have shown that $W$ is a standard \rep of $S_{n+1}$. Our aim now is to find a natural basis for $W$.

Let $a \in W \cap L$ be a non-zero vector, where $L$ is the $(-1)$-eigenspace of $\sigma_{1 (n+1)}$. Note that $\langle a \rangle = W \cap L$. Let us remark here that if we were considering $V_{n-1}$, we would have had $W$ a signed standard representation, and we would have taken $L$ to be the $(+1)$-eigenspace of $\sigma_{1 (n+1)}$ to the same effect.

We write $a = \sum_{i=1}^n a_i$, where $a_i \in E_i$ for each $i$. Now $[\sigma_{1 (n+1)},\sigma]=1$ for each $\sigma \in A_n$ such that $\sigma$ fixes $1$ in the natural action $A_n \curvearrowright \{ 1,2,\ldots,n\}$. Therefore, for each such $\sigma$, \[\sigma(a) \in W \cap L = \langle a \rangle. \] But $W$ is a standard \rep of $A_{n+1}$, and hence $\sigma(a) =a$. So $a_j = \iota_{2j} (a_2)$ for each $j>2$. If $a_1 = \iota_{21} (a_2)$, then in fact $\langle a \rangle \leqslant V_1$ is $A_{n+1}$-invariant, which is a contradiction, since $V_1$ is a sum of standard \reps of $A_{n+1}$. Hence $a_1 \neq \iota_{21} (a_2)$

Let $u = \iota_{21}(a_2) + \sum_{i=2}^n a_i \in V_1$ and set $v_1 = a - u$ and $v_j = \iota_{1j}(a_1) - a_j$ for each $j>1$. Note that $v_i \in E_i$ for each $i$. Now \[ \{ v_1 + u, v_2 + u, \ldots, v_n + u \} \] is a basis for $W$; it is in fact what might be called a natural basis for a standard representation, that is $S_n$ acts by permuting the vectors, and one of the vectors spans $L \cap W$, the $(-1)$-eigenspace of $\sigma_{1 (n+1)}$ in $W$. We can conclude that in particular $\sigma_{1 (n+1)}(v_i + u) = v_i - v_1$ for each $i>1$.

We are going to show that in fact $u=0$. Let us suppose that $u \neq 0$. The strategy now is to find a trivial \rep of $S_{n+1}$ in $V_1$, which will be a contradiction.

We have, for $i>1$,
 \begin{align*} \sigma_{1 (n+1)}(u) &= \sigma_{1 (n+1)}(u + v_i - v_i) \\ &= \sigma_{1 (n+1)}(u + v_i) - \sigma_{1 (n+1)}(v_i) \\ &= v_i -v_1 - \sigma_{1 (n+1)}(v_i) . \end{align*}
But \[ \sigma_{1 (n+1)}(v_i) = \epsilon_1 \rho_{i1} \prod_{j \neq i} \rho_{j1} (v_i) = \epsilon_1 \rho_{i1}(v_i) \in V_0 \oplus E_1 \oplus E_i \] by Lemma~\ref{diamond diagram}. So
\begin{align*} \sigma_{1 (n+1)}(u) &=  v_i -v_1 - \sigma_{1 (n+1)}(v_i) \in V_0 \oplus E_1 \oplus E_i \end{align*}
for each $i \neq 1$. Hence $\sigma_{1 (n+1)}(u) \in V_0 \oplus E_1$. But also $V_1$ is $S_{n+1}$-invariant, and therefore $\sigma_{1 (n+1)}(u) = x_1 \in E_1$. Note that $u \neq 0$ and so $x_1 \neq 0$.

Define $x_i = \iota_{1i}(x_1) \in E_i$ and note that, since $A_{n}$ acts trivially on $u$, $x_i = \sigma_{i (n+1)}(u)$. Now, for $i \neq 1$,
\begin{align*} \sigma_{1 (n+1)}(x_i) &= \sigma_{1 (n+1)} \sigma_{i (n+1)}(u)\\ &= \sigma_{1 (n+1)} \sigma_{i (n+1)} \sigma_{1 i} (u) \\ &=  \sigma_{1 (n+1)} \sigma_{1 i} \sigma_{1 (n+1)} (u) \\ &= \sigma_{i (n+1)}(u) \\ &= x_i. \end{align*}
Note that this calculation is slightly different in the case of $V_{n-1}$ due to extra signs occurring, but the conclusion stays the same.

We have shown that $\{u , x_1 , x_2, \ldots, x_n \}$ forms a basis of a permutation \rep of $S_{n+1}$ within $V_1$. In particular this implies that $\langle u + \sum_{i=1}^n x_i \rangle$ is a one-dimensional \rep of $S_{n+1}$ within $V_1$, which is a contradiction. We conclude that $u=0$.

We have thus shown that a natural basis for $W$ is given by $\{v_1, v_2, \ldots, v_n \}$, and therefore
\begin{align*}
v_2 + v_1 &= \epsilon_1 \sigma_{1 (n+1)}(v_2) \\
&= \prod_{i=2}^n \rho_{i1}(v_2) \\
&= \rho_{21}(v_2). \end{align*}
Also, as $[\epsilon_1 \sigma_{1 (n+1)}, \rho_{21}]=1$,
\begin{align*}
v_3 + v_1 &=\epsilon_1 \sigma_{1 (n+1)}(v_3) \\ &= \epsilon_1 \sigma_{1 (n+1)} \rho_{21} (v_3) \\
&= \rho_{21} \epsilon_1 \sigma_{1 (n+1)}(v_3) \\
&= \rho_{21} ( v_3 + v_1) \\
&= v_3 + \rho_{21}(v_1). \end{align*}
So, combining these two computations with Lemma~\ref{diamond diagram} shows that $\rho_{21} (W) \leqslant V_1$. The same argument works for any $\rho_{ij}$ and any standard \rep $W \leqslant V_1$ of $S_{n+1}$, and these \reps sum up to $V_1$, so we conclude that $\rho_{ij}$ keeps $V_1$ invariant for each $i \neq j$. The same is clearly true for each $\epsilon_i$, and therefore $\Out{n} (V_1) \leqslant V_1$.
Analogously \[\Out{n} (V_{n-1}) \leqslant V_{n-1}. \]

Now we can quotient these two spaces out and obtain a \rep of $\Out{n}$ on the direct sum of $\widetilde{V_0} = V_0/(V_1 \oplus V_{n-1})$ and $\widetilde{V_n} = V_n /(V_1 \oplus V_{n-1})$.

\noindent \textbf{Step 3:} We claim further that $\widetilde{V_0} \oplus \widetilde{V_n}$ is a sum of $\Out{n}$-modules, and the action of $\Out{n}$ on $\widetilde{V_{0}}$ is trivial on $\widetilde{V_n}$ and a sum of determinant representations.

We have shown that $V_0 \oplus V_1$ and $V_{n-1} \oplus V_n$ are $\Out{n}$-invariant, and hence
both $\widetilde{V_{0}}$ and $\widetilde{V_n}$ are representations of $\Out{n}$.
This way we get two maps of the form $\phi  \colon \Out{n} \to \GL_\nu(\K)$, with $\nu \leqslant m$, each of which sends all elements $\epsilon_i$ to either the identity or the minus identity matrix.

Consider the following commutative diagram
\begin{diagram}
\Out{n} & \rTo^{\phi} & \GL_\nu(\K) \\
& \rdTo & \dTo \, \pi \\
& & \mathrm{PGL}_\nu(\K),
\end{diagram}
where $\pi$ is the natural projection. All elements $\epsilon_i$ are in the kernel of the diagonal map, and hence, using Proposition~\ref{press}, we get another commutative diagram
\begin{diagram}
\Out{n} & \rTo & \GL_\nu(\K) \\
\dTo & \rdTo & \dTo \, \pi \\
\GL_n(\Z_2) & \rTo & \mathrm{PGL}_\nu(\K).
\end{diagram}
Now we can use Theorem~\ref{mapfromglnz2}: if $n \geqslant 5$ then the inequality
$\nu \leqslant m < n(n-2) \leqslant 2^{n-1} -1$ allows us to conclude that the bottom map is trivial.
If $n=4$ then we need to additionally use the assumption that $\dim V_1 + \dim V_{n-1} \geqslant 4$. This tells us that \[\nu \leqslant m-4 < n(n-2)-4 \leqslant 2^{n-1} -1\] and hence we can apply the theorem.

In either case, the image of $\Out{n}$ in $\GL_\nu(\K)$ lies in the kernel of $\pi$, which is isomorphic to $\K^*$. So $\phi$ is in fact a sum of identical one-dimensional \reps of $\Out{n}$, and therefore we can apply Lemma~\ref{1-dim reps}. We see that $\phi$ is either a sum of trivial or the determinant representations. But we know the image of $\epsilon_1$ under $\phi$ (depending on whether we are looking at $\widetilde{V_0}$ or $\widetilde{V_n}$), which finishes the proof of this step.

\noindent \textbf{Step 4:} It remains to show that in fact both $V_0$ and $V_n$ are $\Out{n}$-invariant.

Let $v \in V_0$. We know that $\epsilon_1 \sigma_{1 (n+1)}(v) \in V_0$, and therefore in particular its projection onto each of $E_i$ is zero. Now, by Lemma~\ref{diamond diagram}, for $j>1$, the $E_j$ component of
\[ \epsilon_1 \sigma_{1 (n+1)}(v) = \prod_{i=2}^n \rho_{i1} (v) \]
is equal to that of $\rho_{j1}(v)$. Therefore $\rho_{j1}(v) \in V_0 \oplus E_1$ for all $j>1$.

Let $\rho_{21}(v) = v + v'$, where $v' \in E_1$. Hence
\[ \rho_{21}^{-1}(v) = \epsilon_1 \rho_{21} \epsilon_1 (v) = v -v', \]
and so, to ensure that $\rho_{21} \rho_{21}^{-1} = 1$, we need $\rho_{21}(v') = v'$. Now
\[ \epsilon_1 \sigma_{1 (n+1)} (v) = \big( \prod_{i=3}^n \sigma_{i2} \rho_{21} \sigma_{i2} \big) \rho_{21}(v) = v + (n-1)v', \]
which belongs to $V_0$ only if $v'=0$. This shows that $\rho_{21}(v) =v \in V_0$.

The argument works in an identical manner for all $\rho_{ij}$, and for $V_n$. We have therefore finished the proof of this step, and consequently of the proposition.
\end{proof}
\end{prop}

\begin{lem}
\label{beyond n choose 2}
Suppose $\phi  \colon \Out{n} \to \GL(V)$ is an $m$-dimensional $\K$-linear \rep of $\Out{n}$, where $n \geqslant 4$ and $m< {{n+1} \choose 2}$, such that, with the notation of Definition~\ref{reps of W_n}, at least one of $V_2, V_{n-2}$ has non-zero dimension. Then $\phi(\Delta)$ lies in the centre of $\phi(\Out{n})$.
\begin{proof}
\Wlog let us assume that $V_2 \neq \{0\}$.
Lemma~\ref{dimensionofV_i} informs us that \[m - \dim V_2 < {n+1 \choose 2} - {n \choose 2} = n , \] and hence $V_i={0}$ if $i$ is not equal to $0,2$ or $n$.

Now, if $n \geqslant 5$, Lemma~\ref{diamond diagram} shows that each $\rho_{ij}$ preserves $V_0 \oplus V_2$ and $V_n$. Clearly, this is also true for each $\epsilon_i$, and hence $V_0 \oplus V_2$ and $V_n$ are sub\reps of $\Out{n}$. This immediately implies that $\phi(\Delta)$ lies in the centre of $\phi(\Out{n})$, since it lies in the centre of \[\GL(V_0 \oplus V_2) \times \GL(V_n).\]

If $n=4$ then $V = V_0 \oplus V_2 \oplus V_4$, which is precisely the $(+1)$-eigenspace of $\Delta$. Hence, as above, $\phi(\Delta)$ lies in the centre of $\phi(\Out n)$.
\end{proof}
\end{lem}

Combining the two results above yields

\begin{thm}
\label{reps theorem}
\label{thm: B}
Let $\K$ be a field of characteristic equal to zero or greater than $n+1$.
Suppose $\phi  \colon \Out{n} \to \GL(V)$ is an $m$-dimensional $\K$-linear \rep of $\Out{n}$, where $n \geqslant 6$ and $m< {{n+1} \choose 2}$. Then $\phi$ factors through the natural projection $p  \colon \Out{n} \to \GL_n(\Z)$.
\begin{proof}
Firstly, Lemma~\ref{dimensionofV_i} shows that
\[ \forall i \not \in \{1,2,3,n-2,n-1,n\} \colon \dim V_i = 0.\]
We claim that $\phi(\Delta)$ lies in the centre of $\phi(\Out{n})$. We shall consider two cases.

Suppose at least one of $V_2, V_{n-2}$ has non-zero dimension. Then we are in the case of Lemma~\ref{beyond n choose 2}, which asserts the claim.

Suppose now that $V_2 = V_{n-2} = \{0\}$.
Let us note that, since $n \geqslant 6$,
\[ m < { {n+1} \choose2} < n(n-2). \]
We can therefore apply Proposition~\ref{nchoose2case} to $V$ and conclude that, as an $\Out{n}$-module, $V = V_0 \oplus V_1 \oplus V_{n-1} \oplus V_n$. Now $\Delta$ acts as an element of the centre of each $\GL(V_i)$, and hence $\phi(\Delta)$ commutes with $\phi(x)$ for all $x \in \Out{n}$. The claim is thus proven.

The relation $\phi([\Delta,x])=1$ for all $x \in \Out{n}$ in particular holds for $x=\rho_{ij}$, and shows that $\phi(\rho_{ij}) = \phi(\rho_{ij}^\Delta) = \phi(\lambda_{ij})$. Hence we have the following commutative diagram
\[ \begin{diagram}
\Out{n} & \rTo^{\phi} & \GL(V) \\
\dTo^p & \ruTo & \uTo \\
\Out{n}/\langle \! \langle \{ \rho_{ij}\lambda_{ij}^{-1} \mid i \neq j \} \rangle \! \rangle & \cong & \GL_n(\Z)
\end{diagram} \]
which finishes the proof.
\end{proof}
\end{thm}

In a similar vein we obtain

\begin{thm}
Let $\K$ be a field of characteristic equal to zero or greater than 5.
Suppose $\phi  \colon \Out{n} \to \GL(V)$ is an $m$-dimensional $\K$-linear \rep of $\Out{n}$, where $n \in \{4,5\}$ and $m< 2n+1$. Then $\phi$ factors through the natural projection $p  \colon \Out{n} \to \GL_n(\Z)$.
\begin{proof}
First let us suppose that $\dim V_2 + \dim V_{n-2} > 0$. Then we apply Lemma~\ref{beyond n choose 2}, which asserts our claim.

If $V_2 = V_{n-2} = \{0\}$ then either we satisfy the hypothesis of Proposition~\ref{nchoose2case}, in which case we proceed just as in the proof above, or we have $\dim V_1 + \dim V_{n-1} = 2n$. In the latter case, if $n=4$, then $V = V_1 \oplus V_3$ and so
$\phi(\Delta)$ commutes with $\phi(\Out n)$. If $n=5$, then $V = V_1 \oplus V_4$. Lemma~\ref{diamond diagram} tells us that both $V_1$ and $V_4$ are $\Out n$-invariant, and hence in particular $\phi(\Delta)$ lies in the centre of $\phi(\Out n)$.
\end{proof}
\end{thm}

The low-dimensional linear representation theory of $\Out 3$ is the focus of another paper of the author~\cite{kielak2011(2)}.

To put our theorems in context, let us mention the work of Potapchik and Rapinchuk~\cite{potapchikrapinchuk2000}. They study complex linear \reps of $\Aut{n}$ in dimension at most $2n-2$. By using the fact that every \rep of $\Out{n}$ is also a \rep of $\Aut{n}$ via the natural projection $\Aut{n} \to \Out{n}$, we deduce the following statement directly from ~{\cite[Theorem 3.1]{potapchikrapinchuk2000}} of Potapchik--Rapinchuk.

\begin{thm}
Let $\phi  \colon \Out{n} \to \GL_m(\C)$ be a representation, where $n \geqslant 3$ and $m \leqslant 2n-2$. Then $\phi$ factors through the natural projection $p \colon\Out{n} \to \GL_n(\Z)$.
\end{thm}

Theorem~\ref{reps theorem} is a strengthening of the above for large $n$. In the spirit of the work of Potapchik and Rapinchuk we can rephrase it into the following.

\begin{cor}
Let $\phi  \colon \Aut{n} \to \GL_m(\K)$ be a representation over a field $\K$ with characteristic either equal to zero or greater than $n+1$, where $n \geqslant 6$ and $m < {{n+1} \choose 2}$. Then either $\phi$ factors through the natural projection $\Aut{n} \to \GL_n(\Z)$, or it does not vanish on the inner automorphisms of $F_n$.
\begin{proof}
Suppose that $\phi$ does vanish on the inner automorphisms of $F_n$. Then it factors as
\[ \begin{diagram}
\Aut{n} & \rTo^{\phi} & \GL_m(\K) \\
\dTo & \ruTo & \\
\Out{n}
\end{diagram} \]
and the result follows by an application of Theorem~\ref{reps theorem}.
\end{proof}
\end{cor}

\section{Linear representations not factoring through $\Out n \to \GL_n(\Z)$}

In this section we will look into a construction of Grunewald--Lubotzky (see \cite{grunewaldlubotzky2006}) of complex linear representations of dimension $n-1$ of a finite-index subgroups of $\Aut n$ (for $n \geqslant 3$), which we will use to construct representations of $\Out n$ which do not factor through the natural epimorphism \[\Out n \to \GL_n(\Z).\] The only other method of obtaining such representations known to the author is to take the maps \[ \Out n \to \Out m \] constructed by Bridson--Vogtmann~\cite{bridsonvogtmann2011}, and follow them by \[\Out m \to \GL_m(\Z) \to \GL_m(\C).\]

Consider $S$, the set of all epimorphisms $F_n \to \Z_2$, with $F_n = \langle a_1, a_2, \ldots, a_n \rangle$ as before. Note that $\vert S \vert = 2^n - 1$, and that $\Aut n$ acts on $S$. Let $G < \Aut n$ be the stabiliser of \[f  \colon F_n \to \Z_2,\] where
\[ f(a_i) = \left\{ \begin{array}{ccl} 1 & \textrm{ if } & i=n \\ 0 & \textrm{ if } & i \neq n \end{array} \right. .\]
Note that $G$ is of index $2^n - 1$ in $\Aut n$.

\begin{figure}
\begin{center}
\includegraphics[scale=0.5]{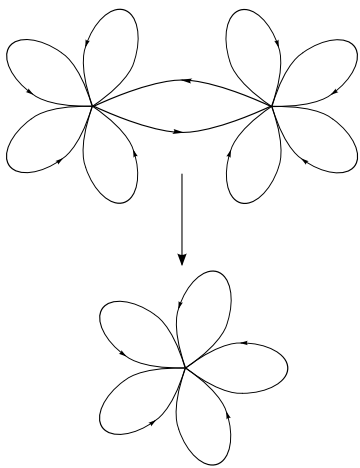}
\caption{The 2-sheeted covering $X \to R_5$}
\label{fig: 2-sheeted covering of the rose}
\end{center}
\end{figure}

Let $R_n$ be the $n$-rose with a fixed isomorphism $\pi_1(R_n) = F_n$, such that the $i^{th}$ petal $b_i$ corresponds to the letter $a_i$. Observe that $G$ contains exactly those based homotopy equivalences of $R_n$ which lift to based homotopy equivalences of a based 2-sheeted covering $X \to R_n$, where $X$ has two vertices joined by lifts of $b_n$, and all the other edges are loops -- see Figure~\ref{fig: 2-sheeted covering of the rose}.

This way we get a map $G \to \Aut {2n-1}$. We can compose it with the natural maps
\[ \Aut {2n-1} \to \GL_{2n-1}(\Z) \to \GL_{2n-1}(\C) \]
to obtain \[\psi \colon G \to \GL_{2n-1}(\C) .\] Since the covering $X \to R_n$ is regular, the action of $G$ on $H_1(X,\C)$ commutes with the action of $\tau$, the non-trivial deck transformation of $X$. Let $V$ denote the $(-1)$-ei\-gen\-space of $\tau$, generated by $\{ \alpha_1, \alpha_2, \ldots, \alpha_{n-1} \}$, where each $\alpha_i$ can be represented by the difference of the two loops in $X$ which project to $b_i$ (see Lemma~\ref{lem: homology of graphs}). We now have
\[ \psi' \colon G \to \GL(V) \cong \GL_{n-1}(\C), \]
which is the representation of Grunewald--Lubotzky.

The group $\Inn n$ of inner automorphisms of $F_n$ is generated by elements \[c_{a_i}  \colon w \to a_i^{-1} w a_i , \] where $w \in F_n$. We immediately see that $\Inn n < G$ and that
\[ \psi'(c_{a_i}) = \left\{ \begin{array}{ccl} \mathrm{I} & \textrm{ if } & i \neq n \\ -\mathrm{I} & \textrm{ if } & i = n \end{array} \right. .\]
We can project $\GL(V) \to \mathrm{GL}(V) / \langle -\mathrm{I} \rangle$ to obtain \[\phi  \colon G/ \Inn n \to \mathrm{GL}(V) / \langle -\mathrm{I} \rangle \cong \mathrm{GL}_{n-1}(\C)/ \langle -\mathrm{I} \rangle . \]
Note that $\vert \Out n: G/ \Inn n \vert = 2^n -1$.

Let $\Torelli_n$ be the kernel of the natural map $\Out n \to \GL( H_1(F_n, \Z))$. It is well known that $\Torelli_n$ is generated by partial conjugations $\rho_{ij} \lambda_{ij}^{-1}$ and commutators $[\rho_{ij},\rho_{ik}]$. We have
\[ \psi(\rho_{ij} \lambda_{ij}^{-1})(\alpha_l) =
\left\{ \begin{array}{ccccl} \alpha_l & \textrm{ if } & j \neq n & \textrm{ or } & l \neq i \\
-\alpha_i & \textrm{ if } & j = n & \textrm{ and } & l = i \end{array} \right. ,\]
and
\[ \psi([\rho_{ij},\rho_{ik}])(\alpha_l) =
\left\{ \begin{array}{ccccl} \alpha_l & \textrm{ if } & n \not\in \{j,k\} & \textrm{ or } & l \neq i \\
\alpha_i - 2 \alpha_k & \textrm{ if } & j = n & \textrm{ and } & l = i  \\
\alpha_i + 2 \alpha_j & \textrm{ if } & k = n & \textrm{ and } & l = i \end{array}
\right. .\]
In particular $\phi(\Torelli_n)$ is infinite.
%

Now to show that this construction has the desired property we will use a standard tool of representation theory, namely Schur functors (see~\cite{fultonharris1991}).

Consider $\mu$, a partition of an even number, and let $\S_\mu$ be the associated Schur's functor. Then $U = \S_\mu V$ is a representation of $\GL_{n-1}(\C)$ factoring through \[\GL_{n-1}(\C) \to \mathrm{GL}(V) / \langle -\mathrm{I} \rangle .\] Thus $U$ is a \rep of $G/ \Inn n$, and we can induce it to a representation \[\theta \colon \Out n \to \GL_m(\C)\] of dimension $m = (2^n - 1) \dim U$.
Note that if $U$ is faithful, then $\Torelli_m \not\leqslant \ker \theta$, and hance $\theta$ does not factor through \[\Out n \to \GL_n(\Z) .\]

When $n \geqslant 4$, the smallest $m$ for which $U$ is faithful is obtained when \[\mu = (1,1) .\] Then $U$ is the second exterior power of $V$, its dimension is ${n-1} \choose 2$, and so \[m = (2^n -1) {{n-1} \choose 2} .\] When $n$ is odd this is smaller then the dimension of the smallest Bridson--Vogtmann representation, and hence the smallest known.

When $n=3$, we need to take $\mu = (2)$, since the second exterior power of $U$ is isomorphic to the determinant \rep in this case. We get \[m =  (2^n -1) {{n} \choose 2} = 21 ,\] which is again smaller that the dimension of the smallest Bridson--Vogtmann representation, which is 55 for $n=3$.

\section{Actions of alternating groups on graphs}
\label{secgraphs}

In this section we establish some lemmata concerning actions of alternating groups on admissible graphs, which will constitute an essential part of our approach.

\begin{dfn}[Admissible graphs]
Let $X$ be a connected graph with no vertices of valence 2, and suppose we have a group $G$ acting on it. We say that $X$ is \emph{$G$-admissible} \iff there is no $G$-invariant non-trivial (i.e. with at least one edge) forest in $X$. We also say that $X$ is \emph{admissible} \iff it is $\mathrm{Aut}(X)$-admissible.
\end{dfn}

Note that in particular a $G$-admissible graph $X$ has no leaves (vertices of valence 1) and no separating edges, independently of what $G$ is.

We prove the following result about admissible graphs:

\begin{lem}
\label{obstruction to admissibility}
Let $X$ be a graph with no separating edges. Suppose $e$ is an edge of $X$ with an endpoint $x$ such that
\[ \forall f \in E(X) \colon e\cap f = \{x\} \Rightarrow m(f) \neq m(e), \]
where $m(f)$ is the minimal length of a simple loop containing an edge $f$. Then $X$ is not admissible.
\begin{proof}
Let $G=\mathrm{Aut}(X)$. Suppose that $X$ is admissible. Then in particular $G . e$ is not a forest. Let $l$ be a simple loop in this orbit. There exists $g \in G$ such that $e \in g .l$. Hence $g.l$ has to contain an edge $f$ of $X$ intersecting $e$ at $x$. But this implies that there exists $h \in G$ such that $h.e = f$. This is a contradiction, since then $m(e) = m(h.e) = m(f)$.
\end{proof}
\end{lem}

The following theorem is due to Marc Culler~\cite{culler1984}, Dmitri Khramtsov~\cite{khramtsov1985} and Bruno Zimmermann~\cite{zimmermann1996} (each independently).

\begin{thm}[(Culler~\cite{culler1984}; Khramtsov~\cite{khramtsov1985}; Zimmermann~\cite{zimmermann1996})]
\label{culler zimmermann}
\label{theorem}
Suppose $G < \Out{m}$ is a finite subgroup. Then there exists a finite $G$-admissible graph $X$ of rank $m$ (with a fixed isomorphism $\pi_1X \cong F_m$), such that $\mathrm{Aut}(X) \leqslant \Out{m}$ contains $G$.
\end{thm}

Since we will be dealing with homology of finite graphs quite frequently in this section, let us observe the following.

\begin{lem}
\label{lem: homology of graphs}
Let $X$ be a finite, oriented graph. Recall that Definition~\ref{graphdef} gives us maps $\iota, \tau  \colon E(X) \to V(X)$. We have the following identification for any field $\K$:
\begin{displaymath}
H_1(X,\K) \cong \Big\{ f  \colon E(X) \to \K \mid \forall a \in V(X)  \colon  \sum_{\iota(x)=a} f(x) = \sum_{\tau(y)=a} f(y) \Big\}.
\end{displaymath}
\begin{proof}
Since $X$ is a CW-complex, we will consider its cellular homology. Since it is 1-dimensional, $H_1(X,\K)$ is the kernel of the boundary map from the $\K$-vector space with basis given by edges of $X$. An element of this vector space is a map $f  \colon E(X) \to \K$, and being in the kernel of the boundary map is equivalent to satisfying the condition \[\sum_{\iota(x)=a} f(x) = \sum_{\tau(y)=a} f(y)\]
at each vertex $a$.
\end{proof}
\end{lem}

We will often refer to each such function $f$ as a choice of \emph{weights} on edges in $X$.

\begin{rmk}
Suppose a group $G$ acts on an oriented graph $X$. Let $f  \colon E(X) \to \K$ be a vector in $H_1(X,\K)$. Then for all $e \in E(X)$
\[ g (f)(e)  = \left\{ \begin{array}{ccc} f(g.e) & \textrm{ if } & o(g.e) = g.o(e) \\ -f(g.e) & \textrm{ if } & o(g.e) \neq g.o(e) \end{array} \right. \]
\end{rmk}

\begin{dfn}
\label{defaxi}
For notational convenience let us define
\begin{displaymath}
\xi= \left\{ \begin{array}{cl} \Delta & \textrm{ if } n \textrm{ is even} \\ \Delta \sigma_{12} & \textrm{ if } n \textrm{ is odd} \end{array} \right.
\end{displaymath}
and $B_n = \langle A_{n+1}, \xi \rangle \leqslant G_n$. We also set $A$ to be either $A_{n-1}$, the pointwise stabiliser of $\{1,2\}$ when $A_{n+1}$ acts on $\{1,2,\ldots,n+1\}$ in the natural way (in the case $n$ is odd), or $A_{n+1}$ (in the case $n$ is even).
\end{dfn}

\begin{lem}
\label{Yisacage}
\label{either rose or cage}
Let $X$ be a connected, oriented, non-trivial graph. Let $n \geqslant 6$. Suppose that $B_n$ acts on $X$ and the action satisfies the following:
\begin{itemize}
\item[(i)] $B_{n}$ acts transitively on the set of (unoriented) edges of $X$;
\item[(ii)] if $A$ acts non-trivially on an edge $e$, then $\xi$ flips each edge in $A.e$ (i.e. it maps the edge to itself, but reverses the orientation);
\item[(iii)] $A$ acts non-trivially on $X$.
\end{itemize}
Then $X$ is either a rose or a cage.
\begin{proof}
Let $e$ be an edge of $X$ such that $A.e \neq e$ as sets (if there was no such $e$, then $A$ would act trivially, since it is perfect). Suppose $e$ is a loop (i.e. is homeomorphic to a circle). Then $X$ is a rose, since it is connected and $B_{n}$ acts on its edges transitively.

Suppose $e$ is not a loop. Suppose further that there exists an edge $f$ which has only one endpoint in common with $e$. Then $f$ cannot be flipped by $\xi$, and in turn must be fixed by $A$, by $(ii)$. This implies that in particular its endpoints are fixed by $A$, and hence also one of the endpoints of $e$ is. Therefore all edges in $A.e$ share a vertex, and, since they all are flipped by $\xi$, they form a cage $C$. Let $\sigma \in A_{n+1}$ be an element taking $e$ to $f$. Then $\sigma$ takes $C$ to a different cage (containing $f$), which is pointwise fixed by $A$, again by $(ii)$. So, $A^\sigma$ has to fix $C$ pointwise. But the intersection $A \cap A^\sigma$ is not empty (since $n>4$), and so the action $A \curvearrowright C$ has a non-trivial kernel. The group $A$ is simple and hence the action has to be trivial. This is a contradiction. We conclude, using the connectedness of $X$, that every edge in $X$ has both endpoints incident with $e$, and therefore $X$ is a cage.
\end{proof}
\end{lem}

In our considerations the following result will be most helpful.

\begin{thm}[({\cite[Theorem 5.2A]{dixonmortimer1996}})]
\label{small index subgroups of A_n}
Suppose $n \geqslant 7$, and let $T < A_n$ be a subgroup of index smaller than ${n+1} \choose 2$. Then $T$ is perfect.
\end{thm}

We are now able to prove the Rose Lemma.

\begin{prop}[(Rose Lemma)]
\label{roseaction}
Suppose $A_{n+1}$ acts on a rose $X$ of rank less than $ {{n+1} \choose 2}$, where $n\geqslant 6$. Then there exists an $A_{n+1}$-invariant choice of orientation of edges of $X$. Moreover, for any field $\K$, the multiplicity of the trivial \rep of $A_{n+1}$ in $V=H_1(X,\K)$ is equal to the number of $A_{n+1}$-orbits of unoriented edges of $X$.
\begin{proof}
Let $e$ be an edge in $X$, and let $T$ be its setwise stabiliser. Then, by the Orbit-Stabiliser Theorem, $\vert A_{n+1}:T \vert < {{n+1} \choose 2}$. Apply Theorem~\ref{small index subgroups of A_n} to $T$ and conclude that it is perfect.

Now, the action of $T$ on $e$ as an oriented edge yields a homomorphism $T \to \Z_2$. Since $T$ is perfect, this homomorphism has to be trivial, and therefore $T$ preserves some (and hence any) orientation of $e$. We can extend this orientation $A_{n+1}$-equivariantly to the orbit of $e$. We can also put weight $1$ on each oriented edge in the orbit, and put weight zero on all other edges of $X$. This way we obtain a non-zero vector $v_e \in H_1(X,\K)$, which is $A_{n+1}$-invariant.

We can repeat the above procedure for each edge in $X$, and conclude the existence of an $A_{n+1}$-invariant orientation on the edges of $X$.

It is clear that $\langle v_f \rangle$ is a trivial \rep of $A_{n+1}$ for each edge $f \in E(X)$. Suppose $v\in H_1(X,\K)$ is a vector spanning a trivial \rep of $A_{n+1}$. Let $e_1, e_2, \ldots e_k$ be a collection of representatives of the edge-orbits of the action of $A_{n+1}$ on $X$. Since $v$ is invariant, it has equal weights on edges in the same orbit. Hence
\[ v = \sum_{i=1}^k v(e_i) v_{e_1} ,\]
where $v(e_i)$ is the weight of $v$ on $e_i$ (with respect to the fixed orientation).
Hence $v \in \langle v_{e_1}, v_{e_2} , \ldots, v_{e_k} \rangle$, and this establishes the equality between number of edge-orbits of $A_{n+1}$ and the multiplicity of the trivial \rep in $H_1(X , \K)$.
\end{proof}
\end{prop}

Let us use similar ideas to prove the following.

\begin{lem}
\label{cagetrivialreps}
Suppose that $A_n$ acts on a non-trivial cage $X$, where $n \geqslant 5$. Then the multiplicity of the trivial representation of $A_n$ in $H_1(X,\C) = V$ is equal to the number of orbits of edges of the cage minus one.

Moreover, if there are at least two edge-orbits, for each edge $e \in E(X)$ we can find an $A_n$-invariant vector $v$ to which $e$ contributes (i.e. the weight of $v$ on $e$ is non-zero).
\begin{proof}
First let us note that $A_n$ has to act on the vertex set of $X$, which gives us a homomorphism $A_n \to \Z_2$. But $A_n$ is perfect, and hence this map has to be trivial. So $A_n$ fixes both vertices of $X$, and therefore preserves the orientation given by choosing one of the vertices to be the image under $\tau$ of all edges.

Suppose that $A_n$ acts transitively on the edges of $X$, and let $v$ be a vector spanning a one-dimensional module in homology. This module has to be a trivial \rep of $A_n$, and so $A_n$ fixes $v$. Therefore $v$ is represented by giving the same weight to each edge. But the sum of weights of outgoing edges has to equal that of ingoing edges at each vertex; in this case it forces the weights to be zero, and therefore $v=0$. This proves our claim in the case when $A_n$ acts transitively on edges of $X$.

Suppose there are at least two orbits of edges in $A_n \curvearrowright X$. Let us label the orbits as $C_0,C_1, \ldots, C_k$. Let us now define vectors $v_i \in H_1(X,\C)$ for $i =1,\ldots, k$ by saying that $v_i$ is represented by giving each edge in $C_i$ weight $\vert C_0 \vert$, each edge in $C_0$ weight $-\vert C_i \vert$, and each edge in $C_j$ weight $0$ for $j \neq 0,i$. Note that each $v_i$ spans a trivial $A_n$-module, and that the vectors $v_i$ are linearly independent. Now let $v$ be a vector in $H_1(X,\C)$ fixed by $A_n$. It necessarily has equal weights on edges in the same orbit; let $\lambda_i$ be the weight of edges in $C_i$. Then we easily verify (using the condition on sums of outgoing and ingoing weights at vertices) that
\begin{displaymath}
\vert C_0 \vert v = \sum_{i=1}^k \lambda_i v_i.
\end{displaymath}
Note that for each edge $e \in E(X)$ there exists an $i$ such that $e \in C_i$ and hence $e$ contributes to $v_i$.
\end{proof}
\end{lem}

\begin{lem}[(Cage Lemma)]
\label{cageaction}
Suppose $A_{n+1}$ (with $n \geqslant 4$) acts on an $m$-cage $X$, so that the action on $V=H_1(X,\C)$ is a sum of standard representations. Assume also that $A_{n+1}$ acts transitively on the edges of $X$. Then in fact $m=n+1$.
\begin{proof}
Let us fix a standard copy of $A_n$ in $A_{n+1}$, i.e. the stabiliser of an element when $A_{n+1}$ acts in a natural way on a set of size $n+1$. We know from the branching rule and our assumption about the representation of $A_{n+1}$, that the multiplicity of the trivial \rep of $A_n$ when acting on $V$ is equal to that of the standard representation.

Suppose that $A_n$ does not fix any edge. Then each orbit gives rise to at least one standard representation of $A_n$. But then, by Lemma~\ref{cagetrivialreps}, we have more standard \reps than trivial representations of $A_n$, which is a contradiction.

Suppose $A_n$ fixes more than one edge. Let $e$ and $e'$ be such edges. Let $\sigma \in A_{n+1}$ be an element sending $e$ to $e'$. Then in particular $\sigma \not\in A_n$ and $A_n ^{\sigma}$ has to fix $e$. Hence $A_{n+1} = \langle A_n, A_n^\sigma \rangle$ fixes $e$, which is a contradiction.

Let $e$ be the unique edge fixed by $A_{n}$, and let $f$ be any other edge of $X$. There exists $\sigma' \in A_{n+1}$ taking $f$ to $e$. So $f$ is the unique fixed edge of $A_{n}^{\sigma'}$, which is a conjugate of $A_n$. We have therefore shown that there is a bijection between edges of $X$ and subgroups in the conjugacy class of $A_n$. There are exactly $n+1$ distinct subgroups of $A_{n+1}$ in the conjugacy class of $A_n$, and hence $m=n+1$.
\end{proof}
\end{lem}

\section{Collapsing maps and the main result}
\label{secmainresult}

In this section we combine the representation theory approach with the graph-theoretic lemmata to prove the main theorem.

\begin{dfn}
\label{collapsingmap}
Let $\pi  \colon X \to X'$ be a surjective morphism of graphs $X$ and $X'$. We say that $\pi$ is a \emph{collapsing map} \iff  for any point $p \in X'$ the preimage $\pi^{-1}(p)$ is connected.
\end{dfn}

Note that this is a slight generalisation of the idea of `collapsing forests', which is present in literature.

\begin{rmk}
Let us observe two facts:
\begin{enumerate}
\item For a graph $X$, giving a subset of $E(X)$ which will be collapsed specifies a collapsing map $\pi$ (up to isomorphism);
\item Any collapsing map $\pi  \colon X \to X'$ induces a surjective map on homology.
\end{enumerate}
\end{rmk}

\begin{dfn}
\label{convenientsplit}
Let $B= A_n \times \Z_2$ for some $n \geqslant 5$, and let $\xi \in B$ denote the element generating the centre of $B$. We say that a \rep $V$ of $B$ \emph{admits a convenient split for $B$} \iff there exists a decomposition $V = U \oplus W$ of $B$-modules, such that, as an $A_n$-module, $U$ is a sum of trivial representations, and such that $\xi$ acts on $W$ as minus identity (the actions of $\xi$ on $U$ and of $A_n$ on $W$ are not prescribed).
\end{dfn}

\begin{figure}
\begin{center}
\includegraphics[scale=1]{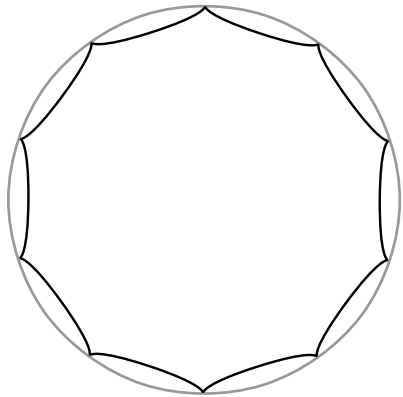}
\caption{A daisy-chain graph -- grey lines represent white edges}
\label{daisy chain}
\end{center}
\end{figure}

\begin{lem}
\label{xiflipsloops}
Let $B= A_n \times \Z_2$ for some $n \geqslant 5$, and let $\xi$ be the generator of the centre of $B$. Suppose that $B$ acts on a graph $X$ so that $A_n < B$ acts non-trivially on each edge of $X$, and such that the action of $B$ on homology admits a convenient split as $H_1(X,\C) = V = U \oplus W$.
Then in fact $\xi$ flips each simple loop in $X$.
\begin{proof}
If $X$ does not contain any simple loops then the result is vacuously true.

Suppose there exists a simple loop $l$ in $X$, and let $v$ be the corresponding vector in homology. We claim that $\xi(v) = -v$, or equivalently that $\xi$ flips $l$.

Suppose for a contradiction that this is not the case. Then $v + \xi (v) \neq 0$, and, as the vector is $\xi$-invariant, it lies in $U$, where $A_n$ acts trivially. So $v + \xi (v)$ is $B$-invariant.

Thus, if $l = \xi.l$ as sets, then $l$ has to be $A_n$ invariant. But $A_n$ cannot act non-trivially on a loop, and hence it fixes each edge. This contradicts our assumption.

Suppose now that we have an edge $f \subseteq l \smallsetminus \xi .l$. In this case we can observe that $A_n .f \subseteq l \cup \xi . l$, since $v + \xi(v)$ is $A_n$-invariant. Note that $A_n . f \subseteq l \cup \xi . l \smallsetminus (l \cap \xi . l)$. Define a collapsing map $\pi  \colon X \to X_f$ by collapsing all edges not contained in the $B$-orbit of $f$. Note that $B$ acts on $X_f$ and $\pi$ is $B$-equivariant. This allows us to use Schur's Lemma (Lemma~\ref{quotientgraphs}) to conclude that $H_1(X_f,\C)$ admits a  convenient split.

We declare the images in $X_f$ of edges of $l$ to be white and images of edges of $\xi.l$ to be black; the action of $\xi$ on $X_f$ will pair up exactly one white edge with exactly one black edge. We claim that $X_f$ has the structure of a daisy-chain graph, where the white edges form a single simple loop, and so do the black edges; see Figure~\ref{daisy chain}.

Let $l'$ be a shortest loop in $X_f$, containing only white edges; we can obtain such a loop since there will be one in the image of $l$. Let $v'$ be the vector corresponding to $l'$ in $H_1(X_f,\C)$. The vector $v' + \xi v'$ is $B$-invariant as before. Moreover, it is not zero, as $v'$ has non-zero weights only on white edges, and $\xi(v')$ has non-zero weights only on black edges. We conclude that $l'$ contains all white edges (since $B$ acts transitively on edges of $X_f$, and $\xi . l$ contains only black edges). We also see that any choice of orientation of $l'$ (i.e. a choice of orientation of its edges such that putting equal weights on each gives a vector in homology) is $B$-invariant; let us fix one such orientation. We can extend it using the action of $\xi$ to a $B$-invariant orientation on the entire graph.

The graph $X_f$ is connected, so there is a vertex of $l'$ from which at least one black edge emanates. But all black edges form a simple loop $\xi . l'$ (since white edges form a simple loop $l'$), and hence in fact we have exactly two black edges emanating from the vertex. The action of $B$ acts transitively on the vertex set of $X_f$ (since it acts transitively on the edge set and preserves the orientation fixed above), so each vertex of $l$ has two white and two black edges emanating from itself. But there are only as many black edges as white, and hence there is a black edge $b$ connecting some two vertices of $l'$. Let $l''$ be a loop formed by $b$ and a shortest subpath of $l'$ connecting the endpoints of $b$; let $v''$ be the corresponding vector in homology. The vector $v'' + \xi v''$ is again $B$-invariant.

Suppose $v'' \neq - \xi v''$. Then, on one hand, $l'' \cup \xi . l''$ contains at most half of all white edges plus one, however on the other hand, being $B$-invariant, it has to contain all white edges. This shows that we have at most two white edges in $l'$, and so at most four edges in $X_f$. But then we would have a non-trivial action of $A_n$, with $n \geqslant 5$, on a set of size 4. This is impossible.

We have shown that $v'' = - \xi v''$, and so in particular $l''$ has length two and contains exactly one white and one black edge. Therefore each black edge in $X_f$ shares both endpoints with a unique white edge. This proves that $X_f$ is a daisy-chain graph as claimed.

Identifying each pair of edges sharing both endpoints gives us an $A_n$-action on a simple loop. Such an action must be trivial, and hence $A_n$ acts on $X_f$ by permuting white and black edges within each pair. This gives us a homomorphism $A_n \to \Z_2^k$ for some $k$. But $A_n$ is perfect, and so each such map must be trivial. Therefore $A_n$ acts trivially on $X_f$. This is a contradiction.

We have therefore shown that $\xi$ sends each simple loop $l$ in $X$ to itself with the opposite orientation. Note that $\C$-linear combinations of simple loops of $X$ span $V$, and so $\xi$ has to act as minus identity on $V$.
\end{proof}
\end{lem}

\begin{lem}
\label{double tree}
Let $X$ be a connected non-trivial graph, on which $\Z_2 = \langle \xi \rangle$ acts in such a way that it flips each simple loop. Then $X = D \cup D'$ as a topological space, where $D$ has the structure of a tree, $D' = \xi . D$, and $D \cap D' = \mathrm{Fix} (\xi)$.
\begin{proof}
Let $F = \mathrm{Fix}(\xi)$ be the fixed point set of $\xi$ in $X$ (where we treat $X$ as a topological space). Let $X' = X \smallsetminus F$.

Firstly, we claim that components of $X'$ are simply connected. Suppose there is a simple loop $l$ in one of the components of $X'$. Since $X' \subseteq X$, $l$ is a simple loop in $X$. As $\xi$ flips all such loops, it flips $l$, and therefore there are two $\xi$-fixed points in $l$. So $l \cap F \neq \emptyset$. This is a contradiction, and therefore each component of $X'$ is simply connected.

We now note that the action of $\xi$ pairs up components of $X'$, and so we can write $X' = \bigsqcup_{i=1}^k (T_i \sqcup T_i ')$ for some $k$, where each $T_i$ is a connected component of $X'$, and $\xi(T_i) = T_i '$. Let $D = \bigsqcup T_i \sqcup F$. Note that $D$ has a structure of a graph: its vertices are vertices of $X$ contained in $D$ together with all points in $F$ which are midpoints of edges in $X$; the edge set is induced by $E(X)$ in an obvious manner.

We now claim that $D$ is in fact a tree. To prove this we will use the following fact: let $p \colon I \to X$ be a path from $x$ to $y$, where $x,y \in D$. We define a path $p'\colon I \to D$ as follows
\[ p'(t) = \left\{ \begin{array}{ccc} p(t) & \textrm{ if } & p(t) \in D \\ \xi . p(t) & \textrm{ if } & p(t) \not\in D \end{array} \right. . \]
Note that $p'$ is a path in $D$ connecting $x$ to $y$. Hence the connectedness of $D$ follows directly from the connectedness of $X$.

Suppose we have a simple loop $l$ in $D$. Then, since $\xi . l = l$ as sets, $l \cap T_i = \emptyset$ for each $i$, and therefore $l \subseteq F$. But then $\xi$ fixes $l$, which contradicts our assumption on $\xi$ flipping all simple loops. So $D$ is a tree.

Define $D' = \bigsqcup T_i' \sqcup F =\xi D$, and note that $F = D \cap D'$ as required.
\end{proof}
\end{lem}

\begin{lem}
\label{double tree 2}
Let $X$ be a connected non-trivial graph, on which $B=A_n \times \Z_2$ acts (with $n \geqslant 5$) in such a way, that there are no $A_n$-fixed edges in $X$. Suppose that $\xi$, the generator of the centre of $B$, flips each simple loop in $X$. Suppose also that all vertices which are not fixed by $A_n$ have valence at least 3. Then in fact all vertices have valence at least $3$, and $A_n$ fixes at most two vertices.
\begin{proof}
Firstly let us apply Lemma~\ref{double tree} to $X$, and conclude that (using notation of the lemma) $X=D \cup D'$. Since $A_n$ commutes with $\xi$, $A_n$ acts on $X/\xi \cong D$. We know that $D$ is a finite tree, and therefore $A_n$ has to fix $d$, its centre (this is a standard fact, see e.g.~\cite{serre2004}). Now let $d'$ be the centre of $D'$. Note that it is possible that $d$ and $d'$ are the same point. Our group $A_n$ acts on $\{ d, d' \}$, and since it is perfect, it has to fix $d$ and $d'$.

Suppose that $A_n$ fixes another point, $x$ say. \Wlog assume that $x \in D$, and take $p$ to be the unique path in $D$ from $x$ to $d$. Now the action of $A_n$ on $p$ can potentially send each subpath of $p$ connecting two points in $F = D \cap D'$ to a subpath lying in $D'$ connecting the same points. Hence the action of $A_n$ on the orbit of $p$ gives a homomorphism $A_n \to \Z_2^k$ for some $k \in \N$. But $A_n$ is perfect, and therefore such a map must be trivial. This implies that $A_n$ fixes $p$, and as $x \neq d$, it has to fix at least one edge. This contradicts our assumption.

We have therefore shown that there are at most two fixed points of $A_n$, namely $d$ and $d'$. Now, if any of these points were of valence less than $3$, then, again as $A_n$ is perfect, each of the edges emanating from it would have to be fixed by $A_n$. This is however impossible, and the proof is finished.
\end{proof}
\end{lem}

We are now ready to prove

\begin{prop}
\label{free reps prop}
Let $n, m \in \N$ be distinct, $n \geqslant 6$, $m< {{n+1} \choose 2}$, and let $\phi \colon  \Out{n} \to \Out{m}$ be a homomorphism. Suppose that the representation $\Out{n} \to \GL(H_1(F_m,\C))=\GL(V)$ induced by $\phi$ satisfies \[V = V_0 \oplus V_1 \oplus V_{n-1} \oplus V_n ,\] with the notation of Definition~\ref{reps of W_n}. Then the image of $\phi$ is contained in a copy of $\Z_2$, the finite group of order two.
\begin{proof}
Before proceeding with the proof,
let us recall some notation, namely Definition~\ref{defaxi}: if $n$ is even, $A=A_{n+1}$ and $\xi = \Delta$; if $n$ is odd, $A=A_{n-1}$ and $\xi = \Delta \sigma_{12}$; we also set $B_n = \langle A_{n+1}, \xi \rangle < G_n$, and $B=\langle A, \xi \rangle \leqslant B_n$.

First let us use Theorem~\ref{theorem} for $\phi(B_n)$ to obtain a finite $B_n$-admissible graph $X$, with an identification $\pi(X) \cong F_m$, such that the action on the conjugacy classes of $F_m$ induced by the action of $B_n$ on $X$ agrees with that given by $\phi$.

The general strategy of this proof will be to first use the results about representation theory of $\Out{n}$ to produce obstructions on the way $B_n$ can act on $X$. Then we will apply the results of this section (dealing with convenient splits), and finally those of Section~\ref{secgraphs}, to conclude that $A < B_n$ has to act trivially on $X$, and hence on the conjugacy classes of elements of its fundamental group. First let us suppose that this last statement is true, and let us deduce the result from there.

\noindent \textbf{Step 0:}
Suppose that $A < B_n$ acts trivially on $X$. We claim that in this case $\phi$ factors through $\Z_2$.

Since $A$ acts trivially on $X$, it acts trivially on the fundamental group of $X$, and hence it lies in the kernel of $\phi$. But $A \leqslant A_{n+1}$, which is simple, and therefore $A_{n+1}$ lies in the kernel of $\phi$. Hence, as $A_{n+1}$ acts transitively on the set $\{ \rho_{ij} \mid i \neq j \}$, $\phi(\rho_{ij}) = \phi(\rho_{jk}),$ and using $[\rho_{ij}^{-1},\rho_{jk}^{-1}]=\rho_{ik}^{-1}$ we see that each $\rho_{ij}$ (and similarly $\lambda_{ij}$) lies in the kernel of $\phi$. This implies that $\phi$ factors through $\Z_2 \cong \langle \epsilon_1 \rangle$.

\noindent \textbf{Step 1:} In what follows let us suppose, for a contradiction, that $A$ does not act trivially on $X$. Let us firstly investigate some of the structure of $X$.

We know that $\Out{n}$ acts via $\phi$ on $H_1(F_m, \C) = V$. Note that also $V \cong H_1(X,\C) \cong \C^m$.

Since $n \geqslant 6$, we have the inequality
\[ m < {{n+1} \choose 2} \leqslant n(n-2),\]
and so we can apply Proposition~\ref{nchoose2case} to $V$.
We claim that $V$ admits a convenient split for $B=\langle A, \xi \rangle$. If $n$ is even, the $A_{n+1}$-modules $V_0 \oplus V_n$ and $V_1 \oplus V_{n-1}$ satisfy the definition.
If $n$ is odd, the sum of all standard representations of $A = A_{n-1}$ is a subspace of $V_1 \oplus V_{n-1}$. Since we chose $A_{n-1}$ to be the stabiliser of 1 an 2 when $A_{n+1}$ acts on $\{1,2,\ldots,n+1\}$, this subspace intersects \[E_1 \oplus E_2 \oplus E_{N \smallsetminus \{1\}} \oplus E_{N \smallsetminus \{1, 2\}}\] trivially. This guarantees that $\xi$ (which equals to $\Delta \sigma_{12}$ in this case) acts on this subspace as minus identity. This proves the claim.

Let us construct a graph $Y$ by collapsing all edges in $X$ which are fixed by $A$ pointwise. Note that, by our assumption, $Y$ is non-trivial (i.e. has at least one edge), and is connected. Since $A$ commutes with $\xi$, we get a $B$-action on $Y$; note that the collapsing map $X \to Y$ is $B$-equivariant. By Lemma~\ref{Schur's lemma}, the $\C$-homology of $Y$ admits a convenient split for $B$. Hence we apply Lemma~\ref{xiflipsloops}, and conclude that $\xi$ flips all simple loops in $Y$.

Note that if we take a vertex $x$ in $Y$ which is not fixed by $A$, then we know that this vertex does not come from collapsing a subgraph of $X$, since we only collapse subgraphs which are $A$-fixed. Therefore such an $x$ comes from a vertex in $X$, and so its valence is at least $3$. This shows that the graph $Y$ (together with the action of $B$ on it) satisfies all conditions of Lemmata~\ref{double tree} and \ref{double tree 2}.

Using the notation of the former lemma, we have $Y = D \cup D'$ and $F=D \cap D' = \mathrm{Fix}(\xi)$. Let $y$ be a point in $\partial D$, i.e. an endpoint of a leaf of $D$. Since its valence (as a vertex of the graph $D \cup D'$) is 2, we see that $y$ is not a vertex of $Y$; it is therefore a midpoint of an edge of $Y$.

Also, the maximal subgraph of $Y$ not containing the $A$-fixed points ($d$ and $d'$) is actually a subgraph of $X$, since any edge collapsed by the map $X \to Y$ yields an $A$-fixed point in $Y$.

\noindent \textbf{Step 2:} We claim that $D$ consists only of leaves. Suppose for a contradiction that this is not the case.

Let $z$ be a farthest (with respect to the graph metric on $D$) vertex of $D$ from $d$, which is not in $\partial D$. We have just assumed that such a vertex is not $d$. Note that $z$ is a vertex of $Y$ and that it cannot be fixed by $A$, as it is neither $d$ nor $d'$. Let $e$ be an edge of $Y$ emanating from $z$, such that its midpoint does not belong to $\partial D$.

Suppose $z \not \in F$. Then all edges in $Y$ emanating from $z$, except for $e$, contain as midpoints points in $\partial D$. There are at least two such edges (since the valence of $z$ is at least 3), and therefore each such edge belongs to a loop of length $2$. Also, neither of these edges forms a loop, since $z \not \in F$, so the shortest loop through any of them is of length 2. This however cannot be true for $e$, since it would require both its endpoints to be in $F$, which is not the case. All of this holds in $X$ as well as $Y$, and we can therefore apply Lemma~\ref{obstruction to admissibility} to $X$ and arrive at a contradiction, since we have assumed that $X$ was $B$-admissible, and hence in particular admissible.

We have thus shown that $z \in F$. But then there exists an edge $f$ in $Y$ emanating from $z$, which is in fact a loop. Note that $f$ is also a loop in $X$. Now consider $X_f$, a graph obtained from $X$ by collapsing all edges but those in $B_n . f$. Note that $B_n$ acts on $X_f$, and the collapsing map $X \to X_f$ is $B_n$-equivariant.

Since $f$ is a loop, $X_f$ is a rose. Also, its rank is at most $m < {{n+1} \choose 2}$. We can therefore apply Proposition~\ref{roseaction} (the Rose Lemma), and obtain an $A_{n+1}$-invariant orientation of edges in $X_f$. By putting equal weight 1 on each edge we obtain an $A_{n+1}$-invariant vector $v \in H_1(X_f,\C)$.

Schur's Lemma (Lemma~\ref{Schur's lemma}) tells us that the image of $V_0 \oplus V_n$ in $H_1(X_f,\C)$ is the sum of all trivial \reps of $A_{n+1}$ in $H_i(X_f,\C)$, and also that the entire group $B_n$ acts trivially on this subspace. Hence $v$ must lie in the image of $V_0 \oplus V_n$, and so $\xi \in B_n$ has to act trivially on it. But $\xi$ flips $f$, which contributes to this vector. This is a contradiction.

This concludes this step, and shows that $D$ is the union of its leaves.

\noindent \textbf{Step 3:} We claim that $X$ is in fact a cage.

We have shown that all edges in $D$ are leaves, and hence are flipped by $\xi$. Hence, in $X$, all edges which are not fixed by $A$ are flipped by $\xi$. Let $f$ be an edge of $X$ flipped by $\xi$, and let $X_f$ be the graph obtained from $X$ by collapsing all edges not contained in $B_n .f$, as before. Note that $A$ acts non-trivially on $f$, since it only fixes one point in $D$.
We can now apply Lemma~\ref{either rose or cage} to $X_f$, which shows that $X_f$ is either a rose or a cage.

The graph $X_f$ cannot be a rose, since if it were, we could construct an $A_{n+1}$-invariant vector $v \in H_1(X_f,\C)$ as in the previous step, on which $\xi$ acts trivially, but to which $f$ (which is flipped by $\xi$) contributes.

So $X_f$ is a cage.
Since $\xi$ flips $f$, it has to permute the two vertices of $X_f$. Also, as $A_{n+1}$ is perfect, it has to fix each of these two vertices.
These vertices have potentially come from non-trivial subgraphs of $X$.
Suppose there exists a simple loop in one of these subgraphs, $l$ say. Let $v$ be a corresponding vector in homology.

Let us assume first that $n$ is odd. We have shown that $\xi$ permutes the vertices of $X_f$ -- in fact this is true for all $\Delta \sigma_{ij}$, since these elements are related by conjugating by elements of $A_{n+1}$. So each $\Delta \sigma_{ij}$ maps $l$ to a loop disjoint from it. So $v + \Delta \sigma_{ij}(v)$ has to be fixed by $A^\sigma$, where $\sigma \in A_{n+1}$ is an element such that $A^\sigma$ commutes with $\Delta \sigma_{ij}$. But each $A^\sigma$ is a simple alternating group, and such groups cannot act on disjoint unions of two circles non-trivially. Hence all $A^\sigma$ fix $l$ pointwise, and therefore so does $A_{n+1}$.

When $n$ is even, $A_{n+1}$ acts trivially on $v + \xi(v)$, and so on a disjoint union of two simple loops $l \cup \xi.l$ as above. So $A_{n+1}$ fixes $l$ pointwise, just as in the odd case.
But then $v \in V_0 \oplus V_n$ and hence the action of $\xi$ on $v$ has to be trivial. This is however not the case.

Therefore the only subgraphs of $X$ we collapsed when constructing $X_f$ were trees. We have however taken $X$ to be $B_n$-admissible, and therefore these trees have to be trivial, i.e. consist of one vertex each.

So $X$ is in fact a cage. Suppose that the action of $A_{n+1}$ on the edge set of $X$ is not transitive. Then, by Lemma~\ref{cagetrivialreps}, there is an $A_{n+1}$-invariant vector $w$ to which $f$ contributes. As $w$ is $A_{n+1}$invariant, it has to lie in the image of $V_0 \oplus V_n$, and hence is $B_n$-invariant. But $\xi$ flips $f$, which is a contradiction. We have thus shown that the action of $A_{n+1}$ on $E(X)$ is transitive. We can apply the Cage Lemma (Lemma~\ref{cageaction}) and conclude that $m=n$, which is a contradiction.
\end{proof}
\end{prop}

The proposition immediately leads to

\begin{thm}
\label{theresult}
\label{thm: A}
Let $n, m \in \N$ be distinct, $n \geqslant 6$, $m< {{n}\choose2}$, and let $\phi \colon  \Out{n} \to \Out{m}$ be a homomorphism. Then the image of $\phi$ is contained in a copy of $\Z_2$, the finite group of order two.
\begin{proof}
As above, let $V = H_1(F_m, \C)$ be a representation of $\Out{n}$ induced by $\phi$. Since $m< {{n}\choose2}$, application of Lemma~\ref{dimensionofV_i} yields
\[ V = V_0 \oplus V_1 \oplus V_{n-1} \oplus V_n ,\]
with the notation of Definition~\ref{reps of W_n} as usual. Hence we can apply Proposition~\ref{free reps prop}, which proves the claim.
\end{proof}
\end{thm}

We can utilise our main tool, Proposition~\ref{free reps prop}, together with a special case of a result of Bridson and Farb~\cite{bridsonfarb2001} to obtain a result reaching a little further. First let us state the required theorem.

\begin{thm}[(Bridson, Farb~\cite{bridsonfarb2001})]
\label{bridson farb}
Suppose $\phi  \colon  \mathrm{PGL}_n(\Z) \to \Out{m}$ is a homomorphism, where $n,m \geqslant 2$. Then the image of $\phi$ is finite.
\end{thm}

Now we can prove

\begin{thm}
\label{main result 2}
\label{thm: C}
Let $n, m \in \N$ be distinct, with $n$ even and at least $6$. Let $\phi \colon  \Out{n} \to \Out{m}$ be a homomorphism. Then the image of $\phi$ is finite, provided that
\[ {n \choose 2} \leqslant m< {{n+1}\choose2}. \]
\begin{proof}
Let $V = H_1(F_m, \C)$ be a representation of $\Out{n}$ induced by $\phi$ as before. Lemma~\ref{dimensionofV_i} shows that either
\[V = V_0 \oplus V_1 \oplus V_{n-1} \oplus V_n \]
or
\[V = V_0 \oplus V_2 \oplus V_{n-2} \oplus V_n. \]
We will proceed by investigating the two cases.

If $V = V_0 \oplus V_1 \oplus V_{n-1} \oplus V_n$, then we can apply Proposition~\ref{free reps prop}, which asserts the claim.

If $V = V_0 \oplus V_2 \oplus V_{n-2} \oplus V_n$, then, as $n$ is even, $\Delta$ acts as identity. Now we can use Theorem~\ref{culler zimmermann} for $\langle \phi(\Delta) \rangle$ and obtain a graph $X$ on which $\Delta$ acts, so that it acts as identity on the homology. Hence in particular it preserves each simple loop (with an orientation) in $X$, since these loops generate the homology. But this implies that $\Delta$ acts trivially on the conjugacy classes in $\pi_1 (X) \cong F_m$, and so $\phi(\Delta) = 1$. This yields the following commutative diagram:
\[ \begin{diagram}
\Out{n} & \rTo^{\phi} & \Out{m} \\
\dTo & \ruTo & \uTo \\
\Out{n}/\langle \! \langle  \Delta \rangle \! \rangle & \cong & \mathrm{PGL}_n(\Z)
\end{diagram} \]
and now an application of Theorem~\ref{bridson farb} finishes the proof.
\end{proof}
\end{thm}

\subsection{Reformulating the statements}

The results contained in this paper can be viewed as an advance in the search for three functions $\alpha, \beta, \gamma  \colon  \N \to \N$, where $\alpha(n)$ is the lowest number such that $\Out{n}$ has an $\alpha(n)$-dimensional complex \rep which does not factor through $\Out{n} \to \GL_n(\Z)$, $\beta(n)$ is the lowest number not equal to $n$ such that there exists a homomorphism $\Out{n} \to \Out{\beta(n)}$ with infinite image (or equivalently a free \rep with infinite image), and $\gamma(n)$ is the lowest number not equal to $n$ such that there exists an embedding $\Out{n} \into \Out{\gamma(n)}$ (or equivalently a faithful free representation).

Our results can then be summarised by saying that if $n \geqslant 6$, then
\[ \alpha(n) \geqslant {{n+1} \choose 2}, \] and
\[ \beta(n), \gamma(n) \geqslant \left\{ \begin{array}{cc} {{n} \choose2} & \textrm{ if $n$ is odd} \\ \\ {{n+1} \choose2} & \textrm{ if $n$ is even} \end{array} \right. . \]

Clearly, $\beta(n) \leqslant \gamma(n)$ for each $n$; it is however unknown if these functions are in fact equal. The relationship between these two functions and $\alpha$ seems to be even more mysterious.

\thispagestyle{plain}
\bibliographystyle{abbrv}
\bibliography{bibliography}

\bigskip

\noindent Dawid Kielak \newline
Mathemathematisches Institut der Universit\"at Bonn \newline
Endenicher Allee 60\newline
D-53115 Bonn\newline
Deutschland \newline
\verb"kielak@math.uni-bonn.de"
\end{document}